\documentclass[11pt]{amsart}

\usepackage{mathpple}

\usepackage{amsmath,amsfonts}
\usepackage{amsrefs}
\usepackage{amsthm}
\usepackage[all]{xy}

\usepackage{hyperref}

\numberwithin{equation}{section}

\newcommand{\C}{\mathbb{C}}

\newcommand{\Gmax}[1]{G_{#1}}
\newcommand{\AW}{E}

\newcommand{\abracket}[1]{\left\langle#1\right\rangle}

\newcommand{\bracket}[1]{\left(#1\right)}

\newcommand{\mc}{\mathcal}

\newcommand{\pa}{\partial}

\renewcommand{\dbar}{\bar\pa}

\newcommand{\A}{\mathcal H}

\DeclareMathOperator{\Res}{Res}

\DeclareMathOperator{\Jac}{Jac}

\theoremstyle{plain}
\newtheorem{thm}{Theorem}[section]
\newtheorem{thm-defn}{Theorem/Definition}[section]

\newtheorem{lem-defn}{Lemma/Definition}[section]

\theoremstyle{definition}

\theoremstyle{remark}

\allowdisplaybreaks[4]  

\begin{document}

 \title{\text{A Mirror Theorem between Landau-Ginzburg models}}
  \author{Si Li}
  \date{}

  \maketitle
{\centering Dedicated to Prof. S.-T. Yau on the occasion of his 65th Birthday \par}


\begin{abstract} We survey on the recent progress toward mirror symmetry between Landau-Ginzburg models. 

\end{abstract}

\section{Introduction}

Calabi-Yau geometry plays a central role for dualities in string theory. The celebrated prediction of counting rational curves on quintic 3-fold \cite{mirror} blows up the mathematical interest on mirror symmetry between symplectic (A-model) and complex (B-model) geometries. A mirror theorem about the equivalence between Gromov-Witten theory and variation of Hodge structures on a large class of mirror Calabi-Yau manifolds has been established by Givental \cite{Givental-mirror} and Lian-Liu-Yau \cite{LLY}.

Calabi-Yau models are closely related to Landau-Ginzburg models \cite{Gepner, vafa-warner, Vafa-LG, Cecotti-LG-CY, chiral-ring, Greene-Vafa-Warner} that are associated to holomorphic functions. Such phenomenon is sometimes called {Calabi-Yau/Landau-Ginzburg correspondence}. The most studied mirror construction involving Landau-Ginzburg models is between toric-varieties and Laurent polynomials \cite{Givental-mirror, Givental-homological, HV}. However, investigation of Landau-Ginzburg mirror pairs was rarely studied in mathematics. This is mainly due to the late appearance of a mathematical theory for Landau-Ginzburg A-models, which was fully established a few years ago by Fan, Jarvis and Ruan \cite{FJR} motivated by work of Witten \cite{Witten}. Now it is commonly referred as FJRW-theory, which can be viewed as a quantum intersection theory on Lefschetz thimbles. See \cite{FT} for a survey. 

The mathematical context of Landau-Ginzburg B-models has a longer history. It dates back to the theory of primitive forms introduced by Saito \cite{Saito-primitive, Saito-residue} around early 1980's as a generalization of the elliptic period integral theory associated to an isolated singularity. Saito's theory leads to systematic examples of Frobenius manifold structure on the universal unfoldings of isolated singularities,  characterizing the genus zero structure of topological Landau-Ginzburg B-model. The analogue for compact Calabi-Yau manifolds is developed in \cite{Barannikov-Kontsevich,Barannikov-quantum}. For higher genus, Givental \cite{G2} proposed a remarkable formula for the total ancestor potential of a semi-simple Frobenius manifold.
The uniqueness of Givental's formula was established by Teleman \cite{T}.
According to Milanov \cite{M} (see also \cite{CI}), Givental's formula can be extended to certain limits of semi-simple locus, which is suffice for our purpose. 

With both sides of Landau-Ginzburg models having firm foundations, it is natural to investigate whether mirror symmetry between Landau-Ginzburg pairs holds. Here the relevant holomorphic functions are weighted homogeneous as required for Landau-Ginzburg A-model. Such polynomials are partially classified by their \emph{central charges}.  Mirror symmetry for ADE singularities (central charge $<1$) is established by Fan-Jarvis-Ruan \cite{FJR} (see also \cite{JKV,FSZ} for A-types). For simple elliptic singularities (central charge $=1$), this is due to Krawitz-Milanov-Shen \cite{KS, MS}. The difficulty to go beyond central charge $=1$ is the lack of computation method for primitive period maps in those cases.  This problem is solved by a recent development of a perturbative theory of primitive forms \cite{LLSaito} (as an analogue of Calabi-Yau cases \cite{Barannikov-quantum,Kevin-Si-BCOV}) together with the help of WDVV equation. Applied to exceptional unimodular singularities (all with $1<$ central charge $<2$), this leads to a mirror theorem \cite{LLSS} for the first nontrivial examples whose central charges exceed $1$. Recently, a general mirror theorem has been proved \cite{H} by a thorough investigation of the method developed in \cite{LLSS}. 

In this note, we will describe a geometric tour toward the general LG/LG mirror theorem based on \cite{H}. This paper is dedicated to Prof. S.-T. Yau on the occasion of his 65th birthday. The author thanks for his sharing of thoughtful views in mathematics and physics and his invaluable support and encouragement all the time. The author would also like to thank Yongbin Ruan for many helpful communications. 

\section{The mirror pairs}
The LG/LG mirror pairs originate from  Berglund-H\"ubsch \cite{BH} that was completed by Krawitz \cite{K}, which is usually called the BHK mirror \cite{CR}. Let $W:\mathbb{C}^N \to\mathbb{C}$ be a weighted homogeneous polynomial with an isolated critical point at the origin. There exist positive rational numbers $q_1, q_2, \ldots, q_N\geq 1/2$ such that
\[
W(\lambda^{q_1}x_1, c^{q_2}x_2, \ldots, \lambda^{q_N}x_N) = \lambda\,W(x_1,x_2, \ldots, x_N), \quad \text{for each} \; \lambda \in \C^\times.
\]
The numbers $q_1, \ldots, q_N$ are called the \emph{weights} of $W$. The \emph{central charge} of $W$, which can be thought of as the ``dimension'' of the LG theory, is defined by
\begin{equation*}
\hat c_W=\sum_{j=1}^{N}(1-2q_j).
\end{equation*}

 We define its \emph{maximal group of diagonal symmetries} to be
\begin{equation*}
\Gmax{W} =
\left\{(\lambda_1,\dots,\lambda_N)\in(\mathbb{C}^{\times})^N\Big\vert\,W(\lambda_1\,x_1,\dots,\lambda_N\,x_N)=W(x_1,\dots,x_N)\right\}.
\end{equation*}
In the BHK mirror construction,  the polynomial $W$ is required to be \emph{invertible} \cite{K, CR}, i.e., the number of variables must equal the number of monomials of $W$. By rescaling the variables, we can always write $W$ as
\begin{equation*}
W=\sum_{i=1}^{N}\prod_{j=1}^{N}x_j^{a_{ij}}.
\end{equation*}
We denote its \emph{exponent matrix} {by $\AW_W =\left(a_{ij}\right)$. The mirror polynomial of $W$ is \cite{BH}
\[
W^T=\sum_{i=1}^{N}\prod_{j=1}^{N}x_j^{a_{ji}},
\]
i.e., the exponent matrix $E_{W^T}$ of the mirror polynomial is the transpose matrix of $\AW_W$. All invertible polynomials have been classified \cite{KreS}: 
a polynomial is invertible if and only if it is a disjoint sum of the three following atomic types, where $a\geq2$ and $a_i\geq2$:
\begin{itemize}
\item \emph{Fermat:}	  $ x^a.$
\item \emph{Chain:}		$x_1^{a_1}x_2 + x_2^{a_2}x_3 + \ldots + x_{N-1}^{a_{N-1}}x_N + x_N^{a_N}.$
\item \emph{Loop:}		$x_1^{a_1}x_2 + x_2^{a_2}x_3 + \ldots + x_N^{a_N}x_1.$
\end{itemize}

Landau-Ginzburg mirror symmetry has to be incorporated with orbifold groups. Let $G\subset \Gmax{W}$ be a subgroup containing 
$
   (e^{2\pi i q_1}, \cdots, e^{2\pi i q_n}). 
$
Then the mirror group $G^T$ for $W^T$ is constructed by \cite{BHe, K}. The BHK mirror pair is 
$$
  (W, G)\leftrightarrow (W^T, G^T). 
$$
Landau-Ginzburg A model of $(W, G)$ is expected to be equivalent to Landau-Ginzburg B-model of $(W^T, G^T)$. When $G=G_W$, $G^T=\{1\}$. Since we do not have a satisfactory theory for orbifold Landau-Ginzburg B-model in general, we will restrict our discussion to $G^T=\{1\}$, i.e. $G=G_W$.

\section{FJRW theory}
FJRW theory \cite{FJR} associates a cohomological field theory (in the sense of \cite{Kont-M}) for a pair $(W, G)$. We consider the case $(W, \Gmax{W})$ with maximal orbifold group. It leads to a \emph{state space} $\A_W$ and a set of linear maps
$$\Lambda_{g,k}^W: (\A_{W})^{\otimes k}\to H^*(\overline{\mathcal M}_{g,k})$$
for  $2g-2+k>0$.
Here $\overline{\mathcal M}_{g,k}$ is the moduli space of stable $k$-pointed curves of genus $g$. The state space is defined as
\[
\A_{W}= \bigoplus_{\gamma \in \Gmax{W}} \A_{\gamma} \qquad \text{where} \qquad \A_{\gamma}:= \left( H^{N_{\gamma}}({\rm Fix}(\gamma), W_{\gamma}^{\infty}; \C )\right)^{\Gmax{W}}.
\]
Here ${\rm Fix}(\gamma)$ is the fixed locus of $\gamma$ and $N_\gamma$ is its dimension as a $\C$-vector space. Furthermore, $W_{\gamma}$ is the restriction of $W$ to ${\rm Fix}(\gamma)$, and $W_{\gamma}^{\infty}$ is ${\rm Re}(W_{\gamma})^{-1}((M, \infty))$ for $M\gg0$. Thus, $\A_{W}$ is dual to the space of Lefschetz thimbles.

The linear maps $\Lambda_{g,k}^W$ are obtained as the moduli space of solutions to the Witten equation
$$
   {\dbar \sigma_i}+\overline{\pa W \over \pa \sigma_i}=0
$$
on Riemann surfaces of genus $g$ with specified boundary conditions via Lefschetz thimbles at $k$ marked points. Here $\sigma_i$'s are sections of suitable orbifold line bundles. See \cite{FJR} for details. Naively,  we can think about $\sigma_i$'s as defining a ``map" (twisted by gravity on the surface $\Sigma_g$) 
$$
   \sigma:  \Sigma_g  ``\to" \C^N
$$ 
satisfying a nonlinear deformation (defined by $W$) of the Cauchy-Riemann equation.  Therefore FJRW theory generalizes the usual Gromov-Witten theory to Landau-Ginzburg models.  See also \cite{CLL} for an algebraic construction in the narrow sectors. 

The \emph{FJRW invariants} are defined by 
$$
    \abracket{\xi_1 \psi_1^{l_1}, \cdots, \xi_k \psi_k^{l_k}}_{g,k}^W=\int_{\overline{\mathcal M}_{g,k}} \Lambda_{g,k}^W(\xi_1, \cdots, \xi_k) \prod_{i=1}^k \psi_i^{l_i}.
$$
Here $\xi_i\in \A_W$, and $\psi_i$ is the $i$-th $\psi$-class on $\overline{\mathcal M}_{g,k}$. The genus zero invariants define a Frobenius manifold structure on $\A_W$ with prepotential function  
$$
\mc F_{0,W}^{FJRW}(\mathbf{t})=\sum_{k\geq 3}{1\over k!} \abracket{\mathbf{t}, \cdots, \mathbf{t}}_{0,k}^W, \quad \mathbf{t}\in \A_W. 
$$

\section{Finding the right B-model}
Let us denote $f=W^T$ in the B-model, with a universal deformation
$$
  F(x,s)=f(x)+ \sum_{\alpha=1}^{\mu} s_\alpha \phi_\alpha(x), 
$$
where $\{\phi_1, \cdots, \phi_\mu\}\in \C[x_1, \cdots, x_N]$ are weighted homogeneous polynomials representing an
additive basis of the Jacobian algebra $Jac(f)$, and $\{s_1, \cdots, s_\mu\}$ parametrizes the deformation
space at the germ $(\C^{\mu},0)$. Saito's primitive form is a family of holomorphic volume forms
$$
\zeta(x,s)=P(x,s)d^Nx, \quad d^Nx=dx_1\wedge\cdots \wedge dx_N,
$$
parametrized by the germ of the deformation. One of the key property of $\zeta(x,s)$ is that it induces a (holomorphic) metric 
$$
   g(\pa_{s_\alpha}, \pa_{s_\beta})=\Res_F\bracket{{\pa_{s_\alpha}F\zeta, {\pa_{s_\beta}F} \zeta}}
$$
which is \emph{flat}. Here $\Res_F$ is the residue pairing associated to the critical points of the holomorphic function $F$ on $\C^N$. Moreover, in terms of the flat affine coordinates
$$
  \tau_\alpha=\tau_\alpha(s), 
$$
the oscillatory integrals satisfy differential equations 
$$
   (\pa_{\tau_\alpha}\pa_{\tau_\beta}-z^{-1}\sum_{\gamma}A_{\alpha\beta}^\gamma(\tau) \pa_{\tau_\gamma}) \int e^{F(x,s(\tau))/z}\zeta(x,s(\tau))=0. 
$$
$A_{\alpha\beta}^\gamma(\tau)$ defines the quantum product for the Frobenius manifold structure on $\Jac(f)$. Together with the flat metric $g$, they determine the potential function of genus zero invariants in the Landau-Ginzburg B-model. 

$\zeta(x,s)$ is constructed by solving a version of Riemann-Hilbert-Birkhoff problem \cite{Saito-primitive}. This abstract nature makes it very difficult to do computations with $\zeta(x,s)$. For weighted homogeneous singularities, the only known expressions are for ADE ($P=1$) and simple elliptic singularities (${1\over P}$ is a period of elliptic curves) \cite{Saito-primitive}. Beyond those, the existence of so-called irrelavant deformations (or negative degree deformations) complicates the situation. 

There is a further subtlety in the B-model.  The primitive forms are not unique. Their moduli space can be identified with the choices of \emph{good basis} \cite{Saito-primitive, Saito-existence, Saito-uniqueness}, whose existence for arbitrary isolated singularity is proved in \cite{Saito-existence}. In the case of weighted homogeneous singularities, a good basis is a homogenous representatives $\{\phi_1, \cdots, \phi_\mu\}\in \C[x_i]$ of an
additive basis of $Jac(f)$, satisfying a tower of higher residue vanishing conditions
$$
  K_f^{(m)}(\phi_\alpha d^Nx, \phi_\beta d^Nx)=0, \quad m=1, 2, \cdots. 
$$
Here $\{K_f^{(m)}\}_{m\geq 0}$ are Saito's higher residue pairing \cite{Saito-residue}, the leading term $K_f^{(0)}$ being the usual residue pairing. It is shown \cite{Saito-primitive} that each choice of a good basis leads to a primitive form, hence a Frobenius manifold structure on $\Jac(f)$. 

Here is a basic example for this phenomenon. Consider a simple elliptic singularity: $f=x_1^3+x_2^3+x_3^3$.  For any $c\in \C$, the following is a good basis 
$$
  \{1, x_1, x_2, x_3, x_1x_2, x_2x_3, x_3x_1, x_1x_2x_3+cf\}. 
$$
In particular, we get a one-parameter family of primitive forms (up to rescaling by a constant).  This is related to the fact that there are two linearly independent periods on elliptic curves. 

There are multiple choices of good basis in the B-model, leading to different genus zero invariants. We have to find the particular one that mirror symmetry favors for. This is identified in \cite{H} as follows. 

\begin{thm}[\cite{H}] Let $f$ be an invertible polynomial of atomic types. Then the following choice $\{\phi_\alpha\}$ is a good basis of $f$. 
\begin{itemize}
\item Let $f=x^a$ be a Fermat, then $\{\phi_{\alpha}\} = \{x^{r} \mid 0\leq r\leq a-2\}$.
\item Let $f = x_1^{a_1}+x_1x_2^{a_2}+\dots+x_{N-1}x_N^{a_N}$ be a chain, then
\[
\{\phi_{\alpha}\} = \left \{\prod_{i=1}^{N}x_i^{r_i}\right\}_{\mathbf{r}} 
\]
 where $\mathbf{r}=(r_1,\cdots,r_N)$ with $r_i\leq  a_i-1$ for all $i$ and $\mathbf{r}$  is not of the form $ (*, \cdots, *, k, a_{N-2l}-1,  \cdots, 0, a_{N-2}-1, 0, a_{N}-1)$ with $k\geq 1$.
\item Let $f=x_1^{a_1}x_N+x_1x_2^{a_2}+\dots+x_{N-1}x_N^{a_N}$ be a loop, then
\[
\{\phi_{\alpha}\} =\left \{\prod_{i=1}^{N}x_i^{r_i} \; \middle | \; 0\leq r_i<a_i \right\}.
\]
\end{itemize}
\end{thm}

Note that we have used the mirror expression for $f=W^T$. For an arbitrary invertible polynomial as a disjoint sum of atomic types, a good basis can be obtained from tensor product of the above good basis from its atomic component. We call this the \emph{standard good basis}. 

The standard good basis appears also in FJRW theory \cite{K}, where Krawitz finds a natural identification of it with elements in FJRW state space $\A_W$. We call this Krawitz's mirror map, which defines a vector space isomorphism between $\A_W$  and $\Jac(W^T)$.  

\section{The mirror theorem}
The works  \cite{T, M} imply that Givental's formula \cite{G2} is also valid for Landau-Ginzburg B-models defined by $W^T$. The full theory of Landau-Ginzburg B-model will be called Saito-Givental theory. Moreover, the higher genus invariants are completely determined by the genus zero data. To establish mirror symmetry between Landau-Ginzburg models, it suffices to show that the Frobenius manifold structures are isomorphic under Krawitz's mirror map. One of the main achievement in \cite{H} is the following reconstruction type theorem. 

\begin{thm}[\cite{H}] Let $W$ be an invertible polynomial with no chain variables of weight $1/2$. Then for both FJRW theory of $(W, G_W)$ and primitive form of $W^T$ with respect to the standard good basis, the genus zero invariants are completely determined by 2-point, 3-point and 4-point functions accompanied with WDVV equation, String equation, Dimension Axiom and Integer Degree Axiom. 
\end{thm}

Here the selection rules of Dimension Axiom and Integer Degree Axiom are natural geometric properties (see \cite{H} for details). The identification of 2-point and 3-point functions is done by Krawitz \cite{K}. This powerful reconstruction theorem therefore reduces the check of mirror symmetry to only 4-point functions (actually a few very special 4-point functions. See \cite{H}). There is a minor situation  for chain types with weight $1/2$ not covered, i.e., $W=x_1^{a_1}x_2 + x_2^{a_2}x_3 + \ldots + x_{N-1}^{a_{N-1}}x_N + x_N^{a_N}$ with $a_N=2$. This is a technical difficulty of missing information about certain FJRW  3-point functions due to the non-algebraic nature of FJRW theory. 

The relevant FJRW 4-point functions can be computed using the method developed in \cite{C, Gu}. Historically, people in the subject of mirror symmetry focused on the computation of 
explicit examples.  One novelty of our computation is that it works for ALL the cases .In the B-model, the difficulty of primitive forms is solved in \cite{LLSaito} by a recursive formula to compute correlation functions up to arbitrary order. These provide enough information for mirror pairs, and it is checked \cite{H} that data from both sides are completely identical! This leads to our mirror theorem. 

\begin{thm}[Landau-Ginzburg Mirror Symmetry Theorem \cite{H}]  Let $W$ be an invertible polynomial with no chain variables of weight $1/2$. Then the FJRW theory of $(W, G_W)$ is equivalent to Saito-Givental theory of $W^T$ at all genera. 

\end{thm}

There are several important questions to explore in this direction. For example, it is desirable to construct the full theory of orbifold Landau-Ginzburg B-model for the pair $(W^T, G^T)$ when $G^T\neq \{1\}$ in order to understand the general mirror pairs. Also, it would be quite interesting to incorporate D-branes for homological mirror symmetry \cite{HMS}. This will be investigated in future works.

\address{\tiny MATHEMATICS SCIENCE CENTER,  TSINGHUA UNIVERSITY, BEIJING, CHINA} \\
\indent \footnotesize{\email{sli@math.tsinghua.edu.cn}}

\end{document}